\documentclass[10pt,reqno,twoside,a4paper]{article}
\usepackage {longtable}
\usepackage{amsmath,amssymb,amsfonts}
\usepackage[T2A]{fontenc}
\usepackage[cp866]{inputenc}
\usepackage[english]{babel}
\usepackage{amsthm}
\usepackage[mathscr]{eucal}
\def\{{\protect\lbrace}
\def\}{\protect\rbrace}

\newcommand{\Hom}{\operatorname{Hom}}

\newcommand{\Ker}{\operatorname{Ker}}

\newcommand{\Sing}{\operatorname{Sing}}

\begin{document}

\begin{center}
\textbf{Modules over Strongly Semiprime Rings}
\end{center}
\begin{center}
A.A. Tuganbaev\footnote{National Research University "MPEI", Lomonosov Moscow State University,\\ e-mail: tuganbaev@gmail.com.}
\end{center}

\textbf{Abstract.} For a ring $A$, the following conditions are equivalent.\\ 
1)~$A$ is a right strongly semiprime ring.\\
2)~Every right $A$-module which is injective with respect to some essential right ideal of the ring $A$, is an injective module;\\
3)~Every quasi-injective right $A$-module which is injective with respect to some essential right ideal of the ring $A$ is an injective module.

The study is supported by Russian Scientific Foundation (project 16-11-10013).

{\bf Key words:} injective module; strongly semiprime ring; quasi-injective module 

\textbf{1. Introduction and preliminaries}

All rings are assumed to be associative and with zero identity element; all modules are unitary and, unless otherwise specified, all modules are right modules. This paper is a continuation of \cite{Tug14}.

For a module $Y$, a module $X$ is said to be \textit{injective with respect to} $Y$ or \textit{$Y$-injective} if for each submodule $Y_1$ of $Y$, every homomorphism $Y_1\to X$ can be extended to a homomorphism $Y\to X$. A module is said to be \textit{injective} if it is injective with respect to each module. A module $X$ is said to be \textit{quasi-injective} if $X$ is injective with respect to $X$. Every finite cyclic group is a quasi-injective non-injective module over the ring of integers $\mathbb{Z}$.

\textbf{Remark 1.1.} The following \textit{Baer's criterion}\footnote{For example, see \cite[Proposition 10.3']{Row88}.} is well-known: if $A$ is a ring and $X$ is a right $A$-module, then $X$ is injective if and only if $X$ is injective with respect to the module $A_A$.

A ring $A$ is said to be \textit{right strongly semiprime} \cite{Han75} if each ideal of $A$ which is an essential right ideal contains a finite subset with zero right annihilator. A ring $A$ is said to be \textit{right strongly prime} \cite{Han75b} if every non-zero ideal of $A$ contains a finite subset with zero right annihilator. It is clear that every right strongly prime ring is right strongly semiprime. The direct product of two finite fields is a finite commutative strongly semiprime ring which is not strongly prime. 

\textbf{Remark 1.2.} Let $A$ be a right strongly prime ring and let $X$ be a right $A$-module. In \cite{Tug14}, it is proved that $X$ is injective if and only if $X$ is injective with respect to some non-zero right ideal of the ring $A$. 

For a module $X$, a submodule $Y$ of $X$ is said to be \textit{essential} in $X$ if $Y\cap Z\ne 0$ for each non-zero submodule $Z$ of $X$. A right module $X$ over the ring $A$ is said to be \textit{non-singular} if the right annihilator $r(x)$ of any non-zero element $x\in X$ is not essential right ideal of the ring $A$. For a module $X$, we denote by $G(X)$ or $\Sing _2X$ the intersection of all submodules $Y$ of the module $X$ such that the factor module $X/Y$ is non-singular. The submodule $G(X)$ is a a fully invariant submodule of $X$ and is called the {\it Goldie radical} or  the \textit{second singular submodule}  of the module $X$.

In connection to Remark 1.1 and Remark 1.2, we prove Theorem 1.3 and Theorem 1.4 which are the main results of this paper.

\textbf{Theorem 1.3.} \textit{For a given ring $A$ with right Goldie radical $G(A_A)$, the following conditions are equivalent.}
\begin{enumerate}
\item[\textbf{1)}] 
\textit{Every non-singular right $A$-module $X$ which is is injective with respect to some essential right ideal of the ring $A$ is an injective module.}
\item[\textbf{2)}] 
\textit{$A/G(A_A)$ is a right strongly semiprime ring.}
\end{enumerate}

\textbf{Theorem 1.4.} \textit{For a given ring $A$, the following conditions are equivalent.}
\begin{enumerate}
\item[\textbf{1)}] 
\textit{$A$ is a right strongly semiprime ring;}
\item[\textbf{2)}] 
\textit{Every right $A$-module which is injective with respect to some essential right ideal of the ring $A$, is an injective module and $A$ is right non-singular.}
\end{enumerate}

\textbf{Remark 1.5.} In connection to Theorem 1.3 and Theorem 1.4, we note that there exist a finite commutative ring $A$, an essential ideal $B$ of the ring $A$, and a non-injective $B$-injective $A$-module $X$. We denote by $A$, $B$ and $X$ the finite commutative ring $\mathbb{Z}/4\mathbb{Z}$, the ideal $2\mathbb{Z}/4\mathbb{Z}$ and the module $B_A$, respectively. Then $B$ is an essential ideal and the module $X$ is injective with respect to $B_A$. Since $X$ is not a direct summand of $A_A$, the module $X$ is not injective.

\textbf{Remark 1.6.} A ring without non-zero nilpotent ideals is called a \textit{semiprime} ring. Every right strongly semiprime ring is a right non-singular semiprime ring \cite{Han75}. The direct product of infinitely many fields is an example of a commutative semiprime non-singular ring which is not strongly semiprime. All finite direct products of rings without zero-dividors and all finite direct products of simple rings are right and left strongly semiprime rings. 

\textbf{Remark 1.7.} A ring $A$ is called a \textit{right Goldie ring} if $A$ is a ring with the maximum condition on right annihilators which does not contain the direct sum of infinitely many non-zero right ideals. If $A$ is a semiprime right Goldie ring, then it is well known \footnote{For example, see \cite[Theorem 3.2.14]{Row88}.} that every essential right ideal of the ring $A$ contains a non-zero-divisor. Therefore, all semiprime right Goldie rings are right strongly semiprime. In particular, all right Noetherian semiprime rings are right strongly semiprime. 

We denote by $\Sing X$ the \textit{singular submodule} of the right $A$-module $X$, that is $\Sing X$ is the fully invariant submodule of $X$ consisting of all elements $x\in X$ such that $r(x)$ is an essential right ideal of the ring $A$. A module $X$ is said to be \textit{singular} if $X=\Sing X$. A module $X$ is called a \textit{Goldie-radical} module if $X=G(X)$. The relation $G(X)=0$ is equivalent to the property that the module $M$ is non-singular. We use well known properties of $\Sing X$ and $G(X)$; e.g., see \cite{Goo76}. A submodule $Y$ of the module $X$ is said to be \textit{closed} in $X$ if $Y=Y'$ for each submodule $Y'$ of $X$ which is an essential extension of the module $Y$.

\textbf{2. The proof of Theorem 1.3 and Theorem 1.4}

The proof of Theorem 1.3 and Theorem 1.4 is decomposed into a series of assertions, some of which are of independent interest. 

\textbf{Lemma 2.1.} \textit{Let $A$ be a ring and let $X$ be a right $A$-module.}
\begin{enumerate}
\item[\textbf{1.}] 
\textit{If $B$ is a right ideal of the ring $A$ and the the modulemodule $X$ is injective with respect to the modulethe the modulemodule $B_A$, then $X$ is injective with respect to the modulethe the modulemodule $(AB)_A$, where $AB$ is the ideal generated by the right ideal $B$. In addition, if the ideal $AB$ contains a finite subset $C$ with $r(C)=0$, then the module $X$ is injective.}
\item[\textbf{2.}] 
\textit{If $X\ne G(X)$, then there exists non-zero right ideal $B$ of the ring $A$ such that the module $B_A$ is isomorphic to a submodule of the module $X$.}
\item[\textbf{3.}]
\textit{Let $Y$ be a right $A$-module, $\{Y_{i}\}_{i\in I}$ be some set of right modules such that the module $X$ is injective with respect to $Y_{i}$ for each $i\in I$ and let $\{f_{i}\in \Hom (Y_{i},Y)\}_{i\in I}$ be some set of homomorphisms. Then the module $X$ is injective with respect to the submodule $\sum _{i\in I}f_{i}(Y_i)$ of the the modulemodule $Y$. In addition, if there exists a monomorphism $A_A\to \sum _{i\in I}f_{i}(Y_i)$, then the module $X$ is injective.} 
\item[\textbf{4.}]
\textit{Let $Y$ be a right $A$-module and let $\{Y_{i}\}_{i\in I}$ be some set of submodules of the module $Y$ such that the module $X$ is injective with respect to $Y_{i}$ for each $i\in I$. Then the module $X$ is injective with respect to the submodule $\sum _{i\in I}Y_i$ of the module $Y$. In addition, if there exists a monomorphism $A_A\to \sum _{i\in I}Y_i$, then the module $X$ is injective.} 
\end{enumerate}

\textbf{Proof.} \textbf{1, 2, 3.} The assertions are proved in Lemma 3, Lemma 4 and Lemma 2.2 from \cite{Tug14}, respectively.

\textbf{4.} The assertion follows from 3 if we denote by $f_i$ the natural embeddings $Y_i\to Y$.~\hfill$\square$

\textbf{Lemma 2.2.} \textit{Let $A$ be a ring, $X$ be a non-singular non-zero right $A$-module, $\{C_i \,|\,i\in I\}$ be the set of all non-zero right ideals of the ring $A$ such that every non-zero submodule of the $A$-module $C_i$ is not isomorphic to a submodule of the module $X$, and let $\{D_j \,|\,j\in J\}$ be the set of all non-zero right ideals $D_j$ of the ring $A$ such that $D_j$ is isomorphic to a submodule of the module $X$. We set $C=\sum _{i\in I}C_i$, $D=\sum _{j\in J}D_j$, and $B=C+D$.}
\begin{enumerate}
\item[\textbf{1.}]
\textit{For any submodule $C'$ of the module $C_A$, every homomorphism $f\colon C'_A\to X$ is the zero homomorphism.} 
\item[\textbf{2.}] 
\textit{The module $X$ is injective with respect to the module $C_A$.}
\item[\textbf{3.}] 
\textit{$B$ is an essential right ideal of the ring $A$.}
\item[\textbf{4.}] 
\textit{If the module $X$ is quasi-injective, then $X$ is injective with respect to the essential right ideal $B$.}
\end{enumerate}

\textbf{Proof.} \textbf{1.} Let us assume that $f\ne 0$. Since $X$ is a non-singular module and $C'/\Ker f\cong f(C')\subseteq X$, we have that $\Ker f$ is not an essential submodule of $C'_A$. There exists a non-zero element $c\in C'$ with $cA\cap \Ker f=0$. Then the non-zero submodule $cA$ of the module $C'$ is isomorphic to the non-zero submodule $f(cA)$ of the module $X$. Therefore, $f(c)\ne 0$. There exists a finite subset $K$ in $I$ such that $c=\sum_{k\in K}c_k$ and $c_k\in C_k$ for all $k\in K$. Since $f(c)\ne 0$, we have that $f(c_k)\ne 0$ for some $k\in K\subseteq I$. Therefore, $c_kA$ is a non-zero submodule of the $A$-module $C_k$ which is isomorphic to a non-zero submodule of the module $X$. This contradicts to the property that $C_k\in \{C_i \,|\,i\in I\}$.

\textbf{2.} The assertion follows from 1.

\textbf{3.} Let us assume that $B$ is not an essential right ideal. Then $B\cap E=0$ for some non-zero right ideal $E$. Then $C\cap E=0$ and $D\cap E=0$. Since $C\cap E=0$, we have that $E\notin \{C_i \,|\,i\in I\}$. Therefore, there exists a non-zero submodule $E_1$ of the module $E$ which is isomorphic to a submodule of the module $X$. Then $E_1\in \{D_j \,|\,j\in J\}$. Therefore, $E_1\subseteq D\cap E=0$. This is a contradiction.

\textbf{4.} Since $X$ is a quasi-injective module, $X$ is injective with respect to any module which is isomorphic to a submodule of the module $X$. Therefore, $X$ is injective with respect to each of the $A$-module $D_j$. By Lemma 2.1(4), the module $X$ is injective with respect to the module $D_A$. In addition, $X$ is injective with respect to the module $C_A$ by 2. By Lemma 2.1(4), the module $X$ is injective with respect to the module $C+D=B$.~\hfill$\square$

\textbf{Proposition 2.3.} \textit{Let $A$ be a right strongly semiprime ring and $X$ be a right $A$-module. If there exists an essential right ideal $B$ of the ring $A$ such that $X$ is injective with respect to the module $B_A$, then $X$ is an injective module.}

\textbf{Proof.} By Lemma 2.1(1), $X$ is injective with respect to the module $(AB)_A$, where $AB$ is the ideal generated by the right ideal $B$. Since $B$ is an essential right ideal and $B\subseteq AB$, the ideal $AB$ is an essential right ideal. Since $A$ is a right strongly semiprime ring, the ideal $AB$ contains a finite subset $K=\{k_1,\ldots ,k_n$ with zero right annihilator $r(K)$. Since $r(K)=r(k_1)\cap \ldots \cap r(k_n)=0$, the module $A_A$ is isomorphic to a submodule of the direct sum of $n$ copies of the module $(AB)_A$. In addition, module $X$ is injective with respect to the module $(AB)_A$. By Lemma 2.1(3), the module $X$ is injective.~\hfill$\square$

For completeness, we briefly prove the following familiar lemma.

\begin{enumerate}
\item[\textbf{1.}]
\textit{If $B$ is an essential right ideal of the ring $A$, then $h(B)$ is an essential right ideal of the ring $h(A)$.} 
\item[\textbf{2.}] 
\textit{If $B$ is a right ideal of the ring $A$ such that $G\subseteq B$ and $h(B)$ is an essential right ideal of the ring $h(A)$, then $B$ is an essential right ideal of the ring $A$.} 
\item[\textbf{3.}] 
\textit{$MG\subseteq G(M)$ for each right $A$-module $M$.}
\item[\textbf{4.}] 
\textit{$XG=0$ and a natural $h(A)$-module $X$ is non-singular. In addition, if $Y$ be an arbitrary non-singular right $A$-module, then $YG=0$ and the $h(A)$-module homomorphisms $Y\to X$ coincide with the $A$-module homomorphisms $Y\to X$. Therefore, $X$ is an $Y$-injective $A$-module if and only if $X$ is an $Y$-injective $h(A)$-module.}
\item[\textbf{5.}] 
\textit{$X$ is an injective $h(A)$-module if and only if $X$ is an injective $A$-module.}
\item[\textbf{6.}] 
\textit{$X_{h(A)}$ is an essential extension of a direct sum of uniform modules if and only if $X_A$ is an essential extension of a direct sum of uniform modules.}
\end{enumerate}

\textbf{Proof.} \textbf{1.} Let us assume that $h(B)$ is not an essential right ideal of the ring $h(A)$. Then there exists a right ideal $C$ of the ring $A$ such that $C$ properly contains $G$ and $h(B)\cap h(C)=h(0)$. Since $h(B)\cap h(C)=h(0)$, we have $B\cap C\subseteq G$. Since $C$ properly contains the closed right ideal $G$, then $C_A$ contains a non-zero submodule $D$ with $D\cap G=0$. Since $B$ is an essential right ideal, $B\cap D\ne 0$; in addition, $(B\cap D)\cap G=0$. Then $h(0)\ne h(B\cap D)\subseteq h(B)\cap h(C)=h(0)$. This is a contradiction.

\textbf{2.} Let us assume that $B$ is not an essential right ideal of the ring $A$. Then $B\cap C=0$ for some non-zero right ideal $C$ of the ring $A$ and 
$G\cap C\subseteq B\cap C=0$. Therefore, $h(C)\ne h(0)$. Since $h(B)$ is an essential right ideal of the ring $h(A)$, we have $h(B)\cap h(C)\ne h(0)$. Let $h(0)\ne h(b)=h(c)\in h(B)\cap h(C)$, where $b\in B$ and $c\in C$. Then $c-b\in G\subseteq B$. Therefore, $c\in B\cap C=0$, whence we have $h(c)=h(0)$. This is a contradiction.

\textbf{3.} For any element $m\in M$, the module $mG_A$ is a Goldie-radical module, since $mG_A$ is a homomorphic image of the Goldie radical module $G$. Therefore, $mG\subseteq G(M)$ and $MG\subseteq G(M)$.

\textbf{4.} By 3, $XG=0$. Let us assume that $x\in X$ and $xh(B)=0$ for some essential right ideal $h(B)$, where $B=h^{-1}(h(B))$ is the complete pre-image of $h(B)$ in the ring $A$. By 2,) $B$ is an essential right ideal of the ring $A$. Then $xB=0$ and $x\in \Sing X=0$. Therefore, $X$ is a non-singular $h(A)$-module.
The remaining part of 4 is directly verified.

\textbf{5.} Let $R$ be one of the rings $A$, $h(A)$ and let $M$ be a right $R$-module. By Lemma 2.1(4), the module $M$ is injective if and only if $M$ is injective with respect to the module $R_R$. Now the assertion follows from 4.

\textbf{6.} The assertion follows from 4.~\hfill$\square$

\textbf{Proposition 2.5.} \textit{Let $A$ be a ring and let $G=A/G(A_A)$. The following conditions are equivalent.}
\begin{enumerate}
\item[\textbf{1)}] 
\textit{Every non-singular right $A$-module $X$ which is is injective with respect to some essential right ideal of the ring $A$ is an injective module.}
\item[\textbf{2)}] 
\textit{Every quasi-injective non-singular right $A$-module $X$ which is is injective with respect to some essential right ideal of the ring $A$ is an injective module.}
\item[\textbf{3)}] 
\textit{Every quasi-injective non-singular right $A$-module is an injective module.}
\item[\textbf{4)}] 
\textit{$A/G(A_A)$ is a right strongly semiprime ring.}
\end{enumerate}

\textbf{Proof.} The implication 1)\,$\Rightarrow$\,2) is obvious.

The implication 2)\,$\Rightarrow$\,3) follows from Lemma 2.2(4).

The equivalence of 3) and 4) is proved in \cite{KutO80}.

4)\,$\Rightarrow$\,1). Let $R$ be one of the rings $A$, $A/G(A_A)$ and let $M$ be a right $R$-module. By Lemma 2.1(4), the module $M$ is injective if and only if $M$ is injective with respect to the module $R_R$.

Let $h\colon A\to A/G$ be a natural ring epimorphism and let $X$ be a non-singular right $A$-module which is injective with respect to some essential right ideal $B$ of the ring $A$. By Lemma 2.4(4), $XG=0$ and $X$ is a natural non-singular $h(A)$-module. By Lemma 2.4(1), $h(B)$ is an essential ideal of the ring $h(A)$. By Lemma 2.4(4), the module $X$ is injective with respect to $h(B)$. By Proposition 2.3, $X$ is an injective $h(A)$-module. By Lemma 2.4(5), $X$ is an injective $A$-module.~\hfill$\square$

\textbf{Remark 2.6. The completion of the proof of Theorem 1.3 and Theorem 1.4.} Theorem 1.3 follows from Proposition 2.5. Theorem 1.4 follows from Proposition 2.3 and the property that every right strongly semiprime ring is right non-singular \cite{Han75}.

\end{document}